%% file: pathplan_arxiv.tex
\documentclass[preprint,12pt,authoryear]{elsarticle}

\newcommand{\Nn}{{\text{N}_\text{n}}}
\newcommand{\Ns}{{\text{N}_\text{s}}}
\newcommand{\Nl}{{\text{N}_\text{l}}}
\newcommand{\Nc}{{\text{N}_\text{c}}}
\newcommand{\Ne}{{\text{N}_\text{e}}}

\usepackage{hyperref}
\usepackage{subfig}
\usepackage{booktabs}
\usepackage{graphicx}
\usepackage[round]{natbib} 
\usepackage{enumitem} 
\usepackage{amssymb}

\journal{Composite Structures}

\include{fabian_macros}

\begin{document}

\begin{frontmatter}

\title{Multi-layer continuous carbon fiber pattern optimization and a spline based path planning interpretation}

\author{Fabian Wein\corref{cor1}\fnref{FAU}}
\ead{fabian.wein@fau.de}
\cortext[cor1]{corresponding author:}
\author{Julian Mirbach\fnref{BA}}
\author{Aniket Angre\fnref{AWI}}
\author{Jannis Greifenstein\fnref{FAU}}
\author{Daniel Hübner\fnref{FAU}}

\affiliation[FAU]{organization={Department of Mathematics, Friedrich-Alexander-Universität Erlangen-Nürnberg}, country={Germany}}
\affiliation[BA]{organization={Broetje-Automation GmbH}, city={Rastede}, country={Germany}}
\affiliation[AWI]{organization={Bionic Lightweight Design and Functional Morphology, Alfred-Wegener-Institut}, city={Bremerhaven}, country={Germany}}

\begin{abstract}

A novel approach for creating tool paths for continuous carbon fiber-reinforced thermoplastic 3D printing is introduced. The aim is to enable load-bearing connections while avoiding non-manufacturable crossings of paths by generating layer specific patterns. We require a graph representation of the structural design with given desired continuous fiber connections. From this, optimal fiber patterns are obtained for each printing layer by solving linear integer optimization problems. Each layer may have a unique solution based on the history of the previous layers. Additionally, a path planning approach is presented which interprets the obtained layers via curves based on quadratic and cubic Bézier splines and their offset curves in constant distance. The single parameter for the construction of the paths is the minimal turning radius of the fibers. The path planning provides a new interpretation for the final geometry of the design to be printed.

\end{abstract}

\end{frontmatter}

\section{Introduction}

Over the past years, 3D printing (3DP)  or additive manufacturing (AM), has gained popularity across various industry sectors ( \cite{VDI:2014:AditiveFertigung} and \cite{AMinAerospace}). It enables the production of parts with complex geometries, see \cite{WZL:Prozesskette}. 3DP is well suited for small batch productions. This makes it especially interesting for applications in aerospace, where batch sizes are typically small.%

Continuous carbon fiber-reinforced thermoplastics (C-CFRTP) are a type of material that consists of thin carbon fibers embedded in a thermoplastic matrix material. This type of material is of particular interest due to its high tensile strength and light weight, see \cite{Wanhill2017}. However, the advantage in tensile strength only holds true for loads carried along the fiber orientation. This anisotropy makes designing and manufacturing parts made from C-CFRTP particularly challenging. Conventional manufacturing processes for C-CFRTP parts, such as \textit{Prepreg} and \textit{Dry-Fiber} only have limited ability to selectively orient fibers (\cite{prepreg}, \cite{DRYFIB}). In this paper we consider an emerging manufacturing process based on the 3D printing of C-CFRTP. This process is called carbon fiber-reinforced printing (CFRP).

At the time of writing, CFRP is still in an early phase towards becoming an established mature process and is subject to active research. In the following we consider CFRP as a subcategory of Fused Deposition Method (FDM), see \cite{fdm} and \cite{fff}, where a filament containing C-CFRTP material is deposited onto a surface by a print head tool. This tool might be attached to an industrial robot.  The printed filament contains thousands of continuous carbon fibers, that are embedded into an isotropic thermoplastic matrix material. Here we denote this as a \textit{fiber bundle}. Before printing, the fiber bundle might have a circular diameter of about 1\,mm, after printing we assume it to have rectangular shape of a width of 2\,mm and a height of 0.2\,mm. Fiber Bundles are printed individually one by one, layer by layer. A brief overview of the full process chain of CFRP with a basic literature selection is given in \cite{Liu:2021:AdditiveComposites}. 

In order to deploy CFRP usefully, tool paths must be defined. Tool paths can be used to program a 3D printing machine. A path defines the way a fiber is placed by the machine. Fibers must be placed such that the desired structural optimized part is created. But the placement must also abide manufacturing constraints. To implement CFRP, the following challenges need to be addressed:
\begin{itemize}
  \item definition of individual C-CFRTP printing paths
  \item constraints on turning radius and parallel alignment of fiber bundles
  \item subsequent accessibility of the printing paths by the print head
  \item possibly necessary support structures 
  \item good mechanical properties of the printed part  
\end{itemize}  
An ideal process would generate a valid printing path respecting these manufacturing constraints in a form, which is transformable to printer machine control commands. Some authors add to the list of tasks the generation of the optimized structure, taking into account the pronounced anisotropy of C-CFRTP. In this work, we contribute to some of the items, but the whole process chain is not yet developed. Our proposed method intends to transform a given structural part into manufacturable paths.

There is still no comprehensive literature in the field of CFRP.  A basic distinguishing feature is the orientation of fiber paths based on the geometry or mechanical load paths. The following are examples of recent publications. In \cite{Papapetrou:2020:StiffnessFramework}, structural optimization based on orthotropic material is performed via a density based and a level-set approach. 
The contribution is the provision of three infill methods for continuous fiber path. The offset method follows the boundary and prefers closed loops. The equally-spaced method creates, despite its name, no parallel path of constant width, which holds also for the third method, a streamline approach. There is no explicit control of fiber curvature. The three infill methods are also numerically evaluated.%

In \cite{Huang:2022:MultiscaleConcurrent}, an anisotropic density based problem formulation is used. The obtained design is interpreted as a graph to generate a single path for the whole structure. Having a single path  prevents closed loops. A curvature formulation is given, but not applied for all numerical examples. Fiber distance is not constant. Interesting is the claimed coupling of the fiber interpretation (via path distance) into the topology optimization, however no details are given. Presumably, the problem is not differentiable and the solution is computationally expensive.

The approach of \cite{Murugan:2022:OrientationField} is rather technical and focuses on the direct generation of machine control code. From a given angle field and domain boundaries, printable path are generated by an explicit Euler based streamline approach. However, minimal path length is not considered. Different layers (stacking) are covered by printing an orthogonal angle field. %

Path planning of anisotropic material is not restricted to fiber-reinforced material. In \cite{Kubalak:2024:OptiToolPath}, anisotropy is assumed for isotropic filament material due to poor layer bonding. The applied  load based streamline approached is described in detail and the numerical results printed. The streamlines are not subject to explicit curvature control.

We believe, that an approach sensitive to the actual loading scenario is better than only interpreting the geometry of a given design, however our approach is significantly different from streamline approaches. The principle considerations for our approach are the following:
\begin{itemize}
  \item we assume that distinct continuous connections are necessary, e.g., from support to load points
  \item not all necessary continuous connections are compatible within a single layer
  \item the user should have the ability to formulate and tune a desired set of continuous connections
\end{itemize} 

\noindent We split our approach into three steps:
\begin{enumerate}[label=(\roman*)]
  \item definition of desired continuous connections in graph form
  \item finding the best realization of (i) via multi-layered fiber pattern optimization
  \item interpretation of the layered fiber patterns (ii) as printable fiber bundle paths    
\end{enumerate}
In this work we assume (i) to be given. In future work we will tune (i) to maximize the strength of the to be printed design interpretation by numerically evaluating (iii). In this work, we use a default configuration of (i) based on the geometry of the given structural design.

To define the continuous connections in (i), a graph representation of the to be printed structure is necessary. It could be obtained by interpreting a given topology optimized design, as it is done in \cite{Huang:2022:MultiscaleConcurrent}. Alternatively, the structural design optimization can even be performed in terms of graph parameters, e.g., via feature-mapping (see the review paper \cite{Wein:2020:Review} and the fiber variants like \cite{Smith:2021:Composite} or \cite{Greifenstein:2023:Spline}). The result of (ii) is an optimal layer wise continuous fiber pattern, but still in a graph setting and not restricted to two dimensions.

Only in the path planning step (iii) do we create a geometric structure that is subject to manufacturing constraints, most importantly strict parallel fiber bundle paths and controlled turning radii. Our proposal for (iii) comes with some restrictions, such as planar two-dimensional layers and no connections with more than four edges, but with some additional effort these restrictions could be relaxed or (iii) could be replaced by a completely different approach. The geometry of the (re)constructed design can be easily controlled by the graph parameters in (i). 

In this paper, we intentionally skip the parametrization of the fiber pattern optimization (i), which allows among other features the encoding of the loading scenarios. This will be presented in forthcoming work.  

This paper is organized as follows: In \secref{sec:int_opt} we present the layered fiber pattern optimization as a linear integer optimization problem. In \secref{sec:pathplanning} the path plan interpretation based on quadratic and cubic Bézier curves is described. A numerical example is shown in \secref{sec:example}. In the final \secref{sec:discussion}, along with the discussion, possible extensions are also formulated. 

\section{Optimal layered fiber patterns}
\label{sec:int_opt}
\begin{figure}[!h]
  \centering
  \includegraphics[width=.75\textwidth]{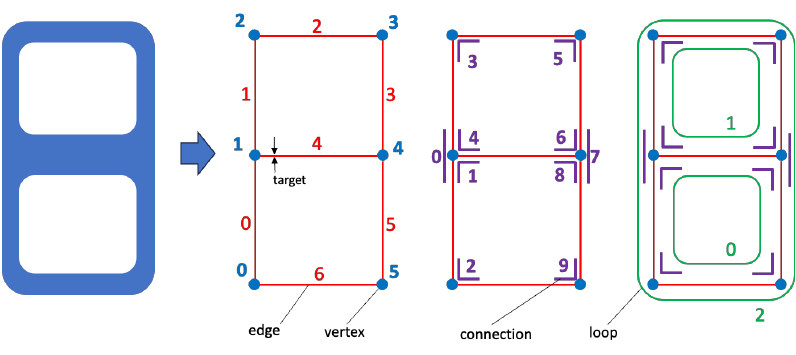} 
  \caption{Setup of the minimalistic example. From a given structure we extract a graph with vertices and edges. Edges form connections and loops are sets of edges/connections.}
  \label{fig:minial_setup}
\end{figure}

We introduce our approach with the minimalistic artificial example in \figref{fig:minial_setup}. The width of the edges is a dimensionless integer representing the number of parallel fiber bundles. The horizontal edge with number 4 has a width of three fiber bundles, all other edges have a width of two fiber bundles. See also the solution in \figref{fig:minimal_three}.

The six vertices define the seven edges $\overline{PQ}$ listed in \tabref{tab:minimal_edges}. The target value represents the number of fiber bundles from the setup.  

\begin{table}[!ht]
\begin{tabular}{lrrrll}
\toprule
id & $P$ & $Q$ & target & (part of connections) & (part of loops) \\
\midrule
0 & 0 & 1 & 2 & 0, 1, 2    & 0, 2 \\ 
1 & 1 & 2 & 2 & 0, 3, 4    & 1, 2 \\
2 & 2 & 3 & 2 & 3, 5       & 1, 2 \\
3 & 3 & 4 & 2 & 5, 6, 7    & 1, 2 \\
4 & 1 & 4 & 3 & 1, 4, 6, 8 & 0, 1 \\
5 & 4 & 5 & 2 & 7, 8, 9    & 0, 2 \\
6 & 0 & 5 & 2 & 2, 9       & 0, 2 \\
\bottomrule
\end{tabular}
\caption{Edges  $\overline{PQ}$ of the setup in \figref{fig:minial_setup}.}
\label{tab:minimal_edges}
\end{table}

We call each unique pair of edges that share a common middle vertex a \textit{connection}. The seven edges correspond to the ten connections given in \tabref{tab:minimal_conns}. In this paper, we set the target value of the connections to the maximum of the targets of both edges. Actually, the target values of the connections are a crucial variable to control the strength of the final solution, but this is beyond the scope of this paper.%

Vertexes and edges are given through the provided graph structure. All connections can be obtained algorithmically from the edges. 

\begin{table}[!ht]
\begin{tabular}{lrrrl}
\toprule
id & edge 1 & edge 2 & target & (part of loops) \\
\midrule
0 & 0 & 1 & 2 & 2 \\
1 & 0 & 4 & 3 & 0 \\
2 & 0 & 6 & 2 & 0, 2 \\
3 & 1 & 2 & 2 & 1, 2 \\
4 & 1 & 4 & 3 & 1 \\
5 & 2 & 3 & 2 & 1, 2 \\
6 & 3 & 4 & 3 & 1 \\
7 & 3 & 5 & 2 & 2 \\
8 & 4 & 5 & 3 & 0 \\
9 & 5 & 6 & 2 & 0, 2 \\
\bottomrule
\end{tabular}
\caption{Connections of the setup in \figref{fig:minial_setup}.}
\label{tab:minimal_conns}
\end{table}

A user needs to provide \textit{loops}, which are ultimately defined by a set of connections. Tabulated in \tabref{tab:minimal_loops}, based on \figref{fig:minial_setup}, loop 0 is the lower small loop, loop 1 is the upper small loop, and loop 2 is the outer large loop. With the intention of printing fiber bundles, we define a loop as the printing unit for layered fiber pattern and also for the successive interpretation through fiber path planning. Whether loops are closed or open has no effect on the presented algorithm. Closed loops do not have a definite beginning or end point. When defining loops via edges rather than connections, it is for the reference implementation necessary to identify a closed loop within the definition by repeating the first edge.

\begin{table}[!ht]
\begin{tabular}{lll}
\toprule
id & connections & edges \\
\midrule
0 & 1, 8, 9, 2 & 0, 4, 5, 6 (,0) \\
1 & 3, 5, 6, 4 & 1, 2, 3, 4 (,1)\\
2 & 0, 3, 5, 7, 9, 2 & 0, 1, 2, 3, 5, 6 (,0)\\
\bottomrule
\end{tabular}
\caption{Loops of the setup in \figref{fig:minial_setup}. When loops are defined by edges, a repeated edge indicates a closed loop for its definition.}
\label{tab:minimal_loops}
\end{table}

The fundamental concept behind optimizing the layered fiber patterns is to consider connections as the crucial parameter but having loops as variables. To achieve this, we describe the relationship between loops and connections as well as loops and edges in the form of matrices called $C \in \mathbb{R}^{\Nc \times \Nl}$ and $E \in \mathbb{N}_0^{\Ne \times \Nl}$, see equation \eqnref{eqn:example_C_cb_E_eb}.  $\Nc$, $\Ne$, and $\Nl$ represent the number of connections, edges, and loops, respectively. Each rows describes for a connection/edge the usage by the loops (three in the example). Connection and edge targets are stored in the vectors $\tilde{c} \in \mathbb{R}^{\Nc}$ and $\tilde{e} \in \mathbb{N}^{\Ne}$. In this work $C$ and $\tilde{c}$ are only integer numbers, but in the general case real numbers allow for more tuning of the optimization problem.
\begin{align}
C = \begin{pmatrix}
  0  &  0  &  1 \\
  1  &  0  &  0 \\
  1  &  0  &  1 \\
  0  &  1  &  1 \\
  0  &  1  &  0 \\
  0  &  1  &  1 \\
  0  &  1  &  0 \\
  0  &  0  &  1 \\
  1  &  0  &  0 \\
  1  &  0  &  1 
\end{pmatrix}\,
\tilde{c} = \begin{pmatrix}
  2 \\
  3 \\
  2 \\
  2 \\
  3 \\
  2 \\
  3 \\
  2 \\
  3 \\
  2
\end{pmatrix}
,\quad
E = \begin{pmatrix}
  1  &  0  &  1 \\
  0  &  1  &  1 \\
  0  &  1  &  1 \\
  0  &  1  &  1 \\
  1  &  1  &  0 \\
  1  &  0  &  1 \\
  1  &  0  &  1
\end{pmatrix}\,
\tilde{e} = \begin{pmatrix}
  2 \\
  2 \\
  2 \\
  2 \\
  3 \\
  2 \\
  2 
\end{pmatrix}.
\label{eqn:example_C_cb_E_eb}
\end{align}

All vectors in this work are column vectors, but when writing a vector with its coefficients in a line we omit the transpose symbol for better readability. %

The connections associated with loop 0 are obtained as
\begin{align*}
  C\,x = (0, 1, 1, 0, 0, 0, 0, 0, 1, 1),\quad x = (1,0,0),
\end{align*}
which gives the first column of $C$. The connections of all loops concurrently are
\begin{align*}
  C\,x = (1, 1, 2, 2, 1, 2, 1, 1, 1, 2),\quad x = (1,1,1),
\end{align*}
what is the sum of all columns of $C$

An obvious problem formulation for a single layer would be to have the actual connections as close as possible to the target vector of connections $\tilde{c}$. This can be written as tracking problem with a discrete loops vector $x \in \mathbb{N}^\Nl_0$ as:

\begin{align*}
  \min_x \; & \| \tilde{c} - C\,x \| \\
       E\,x & \leq \tilde{e}, \\
       x_i & \geq 0\quad \forall \; 1 \leq i \leq \Nl.  
\end{align*}
\noindent Such problems are sometimes written as
\begin{align*}
  \max_x \; & \| C\,x \| \\
       C\,x & \leq \tilde{c}, \\
       E\,x & \leq \tilde{e}, \\
       x_i & \geq 0\quad \forall \; 1 \leq i \leq \Nl  
\end{align*}
Albeit this is an equivalent formulation only when there exists a unique solution $x^*$ with $\| \tilde{c} - C\,x^* \| = 0$, it may serve as a sufficient approximation. The edge constraint is more strict than the constraint for the connections and encompasses it. Thus, the constraint $C\,x \leq \tilde{c}$ can be omitted.

The key concept of layered fiber patterns is the idea of alternating the way fibers are oriented across layers. In order to find an optimal set of fiber patterns, it is necessary to incorporate this concept into the optimization problem. This is done by using the history of optimized loop vectors for previous layers, $x[h] \in \mathbb{N}^{\Nl}_0$ with $h$ being the number of the layer, e.g. $x[h=1]$ describes the fiber pattern of the first layer. This leads to the formulation
\begin{align}
  \min_x \; & \| n\,\tilde{c} - \sum_{h=1}^{n-1} C\,x[h] - C\,x \| \label{eqn:min_tracking} \\
       E\,x & \leq \tilde{e},  \label{eqn:Ex_tracking} \\
       x_i & \geq 0\quad \forall \; 1 \leq i \leq \Nl,  \label{eqn:bounds_tracking}   
\end{align}
with the current layer $n \geq 1$ and $\Nn$ total layers. The three terms in the objective function are the connection target times the layers, $n\,\tilde{c}$, minus the sum of actual connections for the previously solved $n-1$ layers, minus the last term, the connections of the current layer subject to the loops variable, $C\,x$. If we could find in each layer a loop configuration which matches exactly the connection target, the objective term would be zero. As this is practically not the case, underrepresented connections have an increasing growing component in the target which is not sufficiently compensated by the history and current term.

The above optimization problem involves discrete variables, linear constraints and a nonlinear objective function. A possible solver for such problems is Gurobi. However, using a linear program with a scalar product as the objective function is a simpler and more efficient alternative, with readily available open source solvers.   

\subsection{Linear multi-layer optimization}

With a discrete design vector $x \in \mathbb{N}^{\Nl}_0$ for the number of loops in each of the $\Nn$ layers, we formulate the multi-layered fiber pattern problem in it's easiest form as an approximation of the problem \eqnref{eqn:min_tracking} \ldots \eqnref{eqn:bounds_tracking}:
\begin{align}
\text{For }n \in \{1, \ldots, \Nn\} \text{ do:} \nonumber \\ 
\max_{x \in U_\text{ad}} \quad & c(n)^\top x \label{eqn:max_c_x} \\
\text{s.t.} \quad c(n) & = C^\top \left(n \, \tilde{c} - \sum_{h=1}^{n-1} C\,x[h] \right)^p, \label{eqn:single_sheet_c} \\
E \, x & \leq \tilde{e}, \label{eqn:E_x_lt_eb}\\
                x_i & \geq 0 \quad \forall \; 1 \leq i \leq \text{N}_\text{l} \label{eqn:x_gt_0}.
\end{align}
The optimization problem \eqnref{eqn:max_c_x} \ldots \eqnref{eqn:x_gt_0} is to be solved for each layer where the current number of layers is $n$. The admissible design space $U_\text{ad}$ is realized by the constraints \eqnref{eqn:E_x_lt_eb} \ldots \eqnref{eqn:x_gt_0}. For the first layer ($n = 1$), the objective function reduces to the scalar value
\begin{align*}
x^\top  C^\top \tilde{c}^p.
\end{align*}
$c(n)$ is a real vector to size $\Nl$. The exponent $p \in \mathbb{R}^+$ has quite some influence on the optimization problem and is usually chosen as $0 < p \leq 3$. We discuss the impact of $p$ later by example.  

The optimization problems are linear integer problems, which can be solved very efficiently generally to a global optimium. We use \texttt{scipy.optimize.milp} which is based on \cite{higs} implementing \cite{Huangfu:2018:ParallelSimplex}. Technically, we are solving $\min\, -c^\top x$. Note that it is not technically necessary to restrict $C, E, \tilde{c}$ and $\tilde{e}$ to integers. And while only integers make sense for $E$ and $\tilde{e}$, a more sophisticated configuration than our geometric argument can use (positive) real numbers for $\tilde{c}$ and possibly $C$.

A connection represents an uninterrupted bond of edges to be realized by the printing of a continuous fiber enforced fiber bundle in the present printing layer. 
Our motivation is to find a pattern, where single connections are not underrepresented (in view of the previous printing layers).

\begin{figure}[!h]
  \centering
  \includegraphics[width=.45\textwidth]{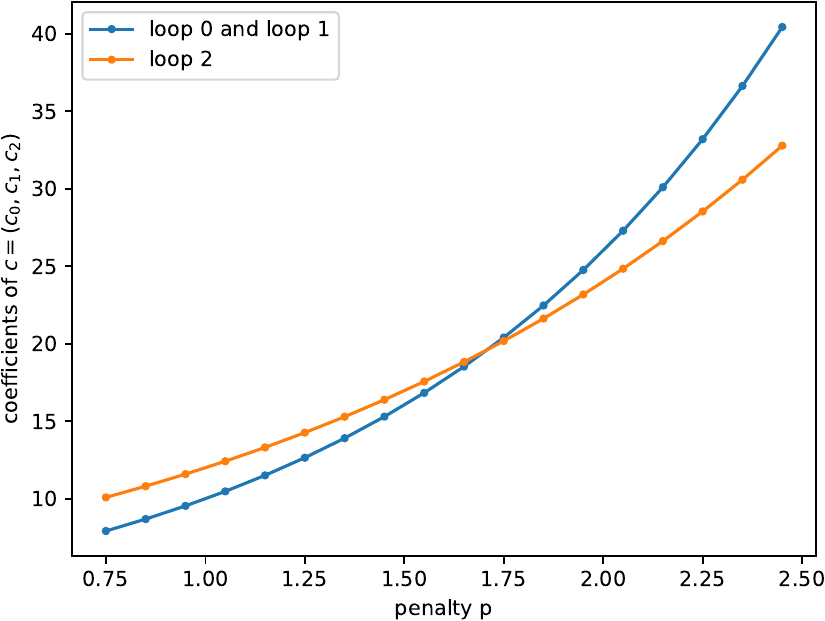}
  \caption{We vary $p$ for $c=C^\top \tilde{c}^p \in \mathbb{R}^{\Nl}$ for the first layer in \figref{fig:minial_setup}. See also \tabref{tab:obj_comb_n_1}.}
  \label{fig:var_power_n_1} 
\end{figure}

\begin{figure}[!h]
  \centering
  \includegraphics[width=.15\textwidth]{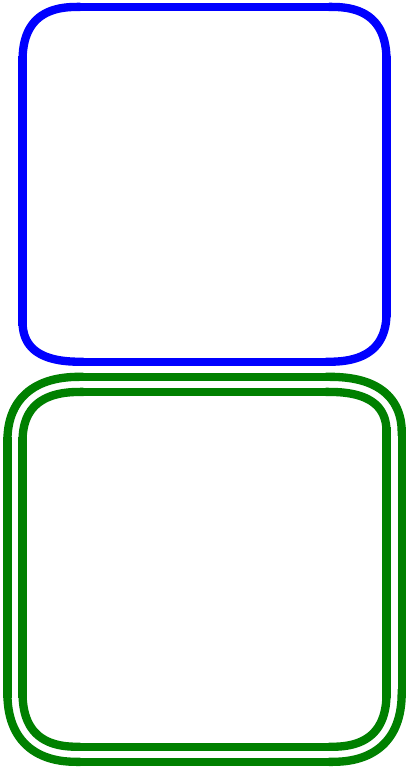} \;
  \includegraphics[width=.15\textwidth]{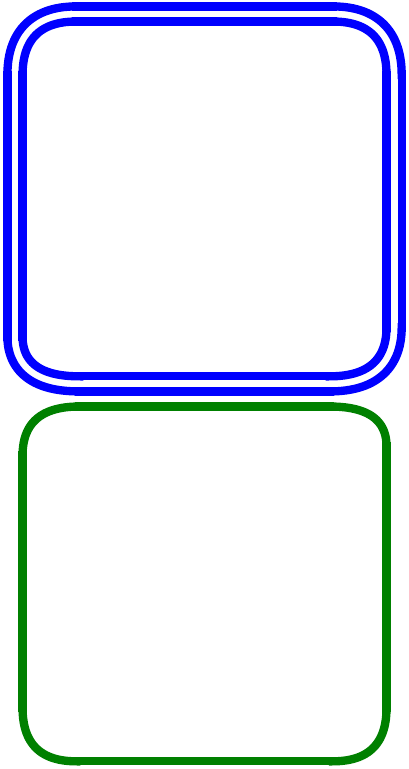} \;
  \includegraphics[width=.15\textwidth]{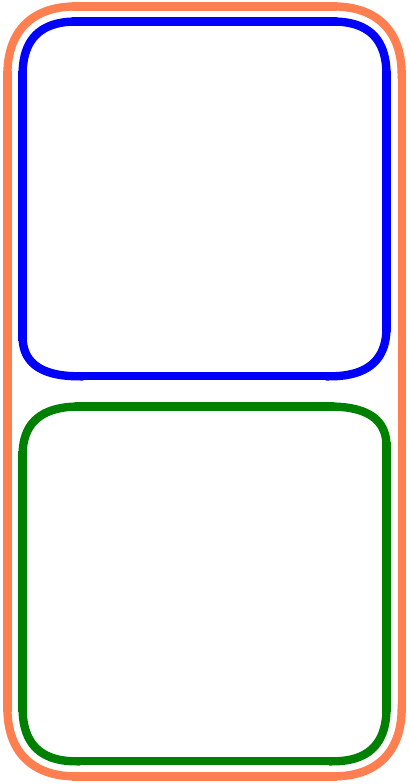}
  \caption{The three relevant layer configurations for the minimal example. For the loops 0, 1 and 2, these are the configurations (2,1,0), (1,2,0) and (1,1,1).}
  \label{fig:minimal_three} 
\end{figure}

For $p=1$ we weight the loops variable $x$ with the vector $C^\top \tilde{c} = (10, 10, 12)$ by counting the target values in $\tilde{c}$ of the connections within these loops. With \tabref{tab:minimal_conns} or \eqnref{eqn:example_C_cb_E_eb} we can confirm that loops 0 and 1 sum up to 10 each and loop 2 amounts to 12. In \tabref{tab:obj_comb_n_1}, we show the objective values for all meaningful loop configurations, e.g. feasible with respect to the constraint \eqnref{eqn:E_x_lt_eb}. See \figref{fig:minimal_three} for the loop configurations and \figref{fig:var_power_n_1} for a variation of $p$ in \eqnref{eqn:single_sheet_c}. 

\begin{table}[!ht]
\begin{tabular}{ccc}
\toprule
$x$ & $C^\top \tilde{c}^1=(10, 10, 12)$ & $C^\top \tilde{c}^2=(26, 26, 24)$ \\
\midrule
(2,1,0) & 30 & 78 \\
(1,2,0) & 30 & 78 \\
(1,1,1) & 32 & 76 \\
(0,0,2) & 24 & 48 \\
\bottomrule
\end{tabular}
\caption{Objective function ($x^\top c$) values for possible loop combinations at the first layer for $c=C^\top \tilde{c}^1=(10, 10, 12)$ and $c=C^\top \tilde{c}^2=(26, 26, 24)$.}
\label{tab:obj_comb_n_1}
\end{table}

We note that the larger loop 2 possesses two additional connections compared to loops 0 and 1 (refer to \tabref{tab:minimal_loops}). However, both smaller loops share the horizontal edge 4 with target three, thereby contributing two connections with target three for loop 0 and 1. The doubling of these smaller loops in the (2,1,0) and (1,2,0) configuration leads to the results portrayed in \tabref{tab:obj_comb_n_1} and \figref{fig:minimal_three}. If the targets of the connections were determined by the minimum of the targets of the edges or if the target for edge 4 were two, the configuration $(1,1,1)$ would outperform any other configuration for any value of $p$. %

The optimal six-layer configurations for the minimal example with $p=2$ are presented in \tabref{tab:x_6_layers}.

As discussed above (\tabref{tab:obj_comb_n_1}), a small loop shows for the initial layer the best performance and the optimal solution consists of twice loop 0, once loop 1 and no large loop 2. 
For the second layer, the two connections with a weight of three for loop 0 contribute to $(2\,3 - 2)^2$ in \eqnref{eqn:single_sheet_c}, while the weight for loop 1 is $(2\,3 - 1)^2$.
The optimal arrangement is (1,2,0), which mirrors the first layer. For the third layer, using solution (1,1,1), the large loop 1 incorporates the vertical connections 0 and 7. This arrangement is repeated for the second set of three layers. \tabref{tab:final_c_6_layers} displays the accumulated sum of selected connections for all six layers. The vertical connections 0 and 7 are only utilized with a single loop in layers 3 and 6, resulting in a sum of two as opposed to $6 \cdot 2=12$. This may appear counterintuitive and unbalanced; however, it results in a minimalistic example with an edge target of only two. A smaller value of $p$ significantly affects the solution, which is a permutation of the configurations (2,1,0), (1,2,0), and (1,1,1).

\begin{table}[!ht]
\begin{tabular}{rrr}
\toprule
$n$ & $x^*$  & $c$ \\
\midrule
1 & (2,1,0) & (26,26,24) \\ 
2 & (1,2,0) & (40,48,58) \\
3 & (1,1,1) & (90,90,108) \\
4 & (2,1,0) & (146,146,134) \\ 
5 & (1,2,0) & (180,232,212) \\
6 & (1,1,1) & (274,274,306) \\
\bottomrule
\end{tabular}
\caption{Optimization of the minimal example for the first six layers with $p=2$.}
\label{tab:x_6_layers}
\end{table}

\begin{table}[!ht]
\begin{tabular}{rrrrrrr}
\toprule
cid & $A$ & $C$ & $B$ & loops & $\Sigma$ vs. $n\,\tilde{c}$ & usage in layers \\
\midrule
  0 & $(   0,   0)$&$(   0, 100)$&$(   0, 200)$ & 2 &  2 vs.   12 & 0 0 1 0 0 1 \\ 
  1 & $(   0,   0)$&$(   0, 100)$&$( 100, 100)$ & 0 &   8 vs.   18 & 2 1 1 2 1 1 \\ 
  2 & $(   0, 100)$&$(   0,   0)$&$( 100,   0)$ & 0, 2 &  10 vs.   12 & 2 1 2 2 1 2 \\ 
  3 & $(   0, 100)$&$(   0, 200)$&$( 100, 200)$ & 1, 2 &  10 vs.   12 & 1 2 2 1 2 2 \\ 
  4 & $(   0, 200)$&$(   0, 100)$&$( 100, 100)$ & 1 &   8 vs.   18 & 1 2 1 1 2 1 \\ 
  5 & $(   0, 200)$&$( 100, 200)$&$( 100, 100)$ & 1, 2 &  10 vs.   12 & 1 2 2 1 2 2 \\ 
  6 & $( 100, 200)$&$( 100, 100)$&$(   0, 100)$ & 1 &   8 vs.   18 & 1 2 1 1 2 1 \\ 
  7 & $( 100, 200)$&$( 100, 100)$&$( 100,   0)$ & 2 &   2 vs.   12 & 0 0 1 0 0 1 \\ 
  8 & $(   0, 100)$&$( 100, 100)$&$( 100,   0)$ & 0 &   8 vs.   18 & 2 1 1 2 1 1 \\ 
  9 & $( 100, 100)$&$( 100,   0)$&$(   0,   0)$ & 0, 2 &  10 vs.   12 & 2 1 2 2 1 2 \\ 
\bottomrule
\end{tabular}
\caption{Final sum of connections for optimization of six layers with $p=2$, see \tabref{tab:x_6_layers}. For convenience the connections are given by three vertices $A-C-B$.}
\label{tab:final_c_6_layers}
\end{table}

\subsection{Multiple sheets}
\begin{figure}[!h]
  \centering
  \subfloat[junction]{
  \includegraphics[width=.2\textwidth]{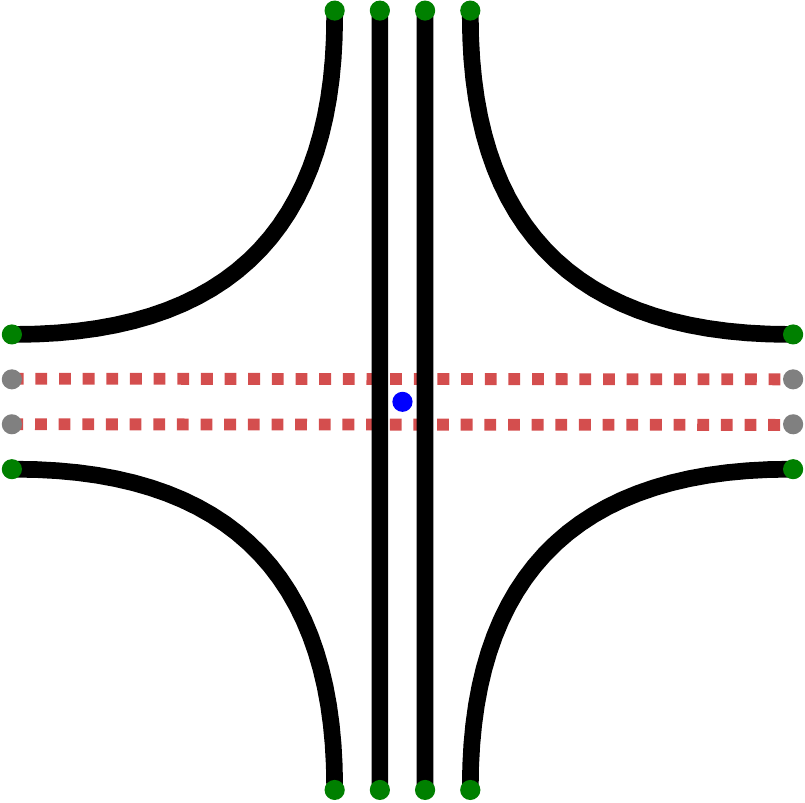}
  \label{fig:crossing_4} }\,
  \subfloat[two sheets to separate incompatible loops]{
  \includegraphics[width=.35\textwidth]{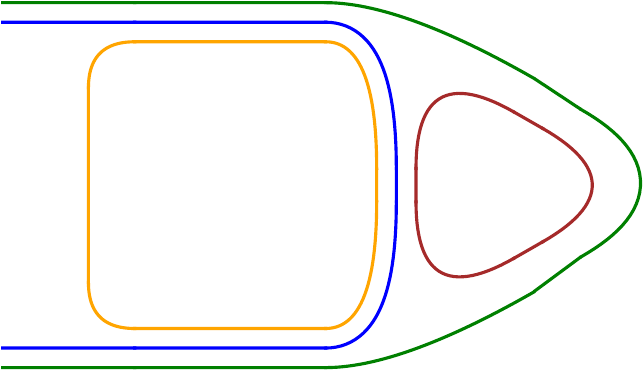} \,
  \includegraphics[width=.35\textwidth]{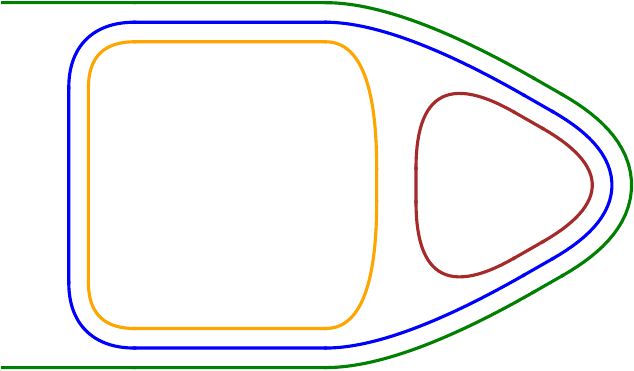}
  \label{fig:two_sheets} }
  \caption{For the junction of four edges, the dashed red horizontal crossing is not concurrently possible with the vertical crossing. For a variant of the minmalistic example with support extensions, two sheets sheets are necessary, where the blue loop is unique in each sheet, the other loops are repeated.}
  
\end{figure}

We define junctions as vertices, where two or more edges meet. A junction with only two edges is considered trivial. For junctions with three edges, all three potential connections are always concurrently possible. However, for junctions with four or more edges, crossing connections are created by non-neighboring edges. A four-edge junction has two different directions of crossing connections, but only one crossing direction is possible at a time in a single layer, as shown in \figref{fig:crossing_4}. Similar is the situation with overlapping loops, see \figref{fig:two_sheets} for an example.

As remedy we introduce \textit{sheets}, each of which contains a configuration of compatible loops. In an abstract way, a sheet represents what is called an admissible set in optimization. The loop relationship matrices become sheet dependent as $C^s$ and $E^s$ with the sheet index $s$ and $\Ns$ the number of sheets. One can either define unique loop indices with constant matrices but sheet dependent loop bounds  or have sheet dependent loop definitions where each sheet has $\Nl^s$ loops but edges and connections are globally indexed. We realized the second approach, but the other approach also fits in the given formula framework.

For each layer $n$ we search the best solution over each sheet $s$. To this end we maximize in each layer each sheet and store the index of the best one. 
\begin{align}
\text{For }n \in \{1, \ldots, \Nn\} \text{ do:} \nonumber \\ 
\max_{s \in \{1, \ldots, \Ns\}} \quad \max_{x \in U_{\text{ad}_s}} \quad & c(s,n)^\top x \label{eqn:max_max_cTx} \\
\text{s.t.} \quad  c(s,n)& = {C^s}^\top \left(n \, \tilde{c} - \sum_{h=1}^{n-1} C^{s[h]}\,x[h] \right)^p,\label{eqn:multi_sheet_c} \\
          E^s \, x & \leq \tilde{e}, \label{eqn:Es_x_lt_eb} \\
                 C^s \, x & \geq \underline{c}, \label{eqn:Cs_x_gt_cb}\\
                 V^s \, x & \leq \overline{v}, \label{eqn:Vs_x_lt_vb}\\
                 V^s \, x & \geq \underline{v}, \label{eqn:Vs_x_gt_vb}\\
                      x_i & \geq 0\quad \forall \; 1 \leq i \leq \Nl^s. \label{eqn:xs_bounds}
\end{align}
$s[h]$ is the stored sheet index for the history layer $h$. The admissible space $U_{\text{ad}_s}$ is now dependent on the sheet and defined by the constraints \eqnref{eqn:Es_x_lt_eb} \ldots \eqnref{eqn:xs_bounds}. Note, that the number of layer $\Nl^s$ is also sheet dependent, hence also the size of $x$ and $c(s,n)$. The constraints \eqnref{eqn:Cs_x_gt_cb} \ldots \eqnref{eqn:Vs_x_gt_vb} are optional. An application of the minimal connections bound \eqnref{eqn:Cs_x_gt_cb} is to have a minimal connection for all connected neighboring edges, as depicted in \figref{fig:crossing_4}. 
The assumption is that in the printing process, the required filling of matrix material within junctions is easier to print with a fiber bundle at the boundary. The vertex constraints \eqnref{eqn:Vs_x_lt_vb} and  \eqnref{eqn:Vs_x_gt_vb} have less obvious usage and are merely given for completeness. They can control the minimal and maximal number of edges entering to a junction via $\underline{v}$ and $\overline{v}$ respectively. 

Already the choice of loops for the single sheet case has a crucial effect on the obtained optimization result. Same is true for the multiple-sheets configuration. When a certain loop, e.g., one with a specific crossing direction for a junction with at least four edges, is only contained in one sheet, the other loops of this sheet need to be sufficiently ``attractive'', too, such that the sheet will represent the optimal solution over all sheets. Hence, the repetition of compatible loops over several sheets is necessary. It is a common case, that a sheet differs only for a single loop from other sheets. The identification of loops is not within the scope of this paper. It might be based on combinatorial aspects (e.g. mathematical power sets) or on load paths from finite element analysis. For both cases, a feedback from numerical validation of the interpreted path plan appears essential.

Alternatively to sheets, incompatible connections within a layer can be prevented by nonlinear constraints.

The fiber pattern optimization is not restricted to two dimensions and loops are allowed to connect arbitrary vertices in three dimensional space, however the interpretation and manufacturing approach is beyond the scope of this paper. Infeasible configurations need to be avoided by proper sheet definitions.

\section{Fiber path planning} 
\label{sec:pathplanning}
In this section, we outline our approach for interpreting the optimization result of the optimized layered fiber pattern by parallel curves subject to a controlled radius. With additional processing, this allows the generation of printer control commands to print the fiber bundles. However, in the context of this study, our focus is solely on visualizing the numerical results from the fiber pattern optimization. Each layer is interpreted individually.

Unlike fiber pattern optimization, this path planning algorithm is limited to two-dimensional layers. One possible extension would be to transform the layers onto a curved surface. 

Generally, we assume all edges are filled with fibers up to the target number. This requirement cannot be fulfilled with for most junctions of three and more edges, therefore filling empty loop space by fiber bundles or filling material is necessary. Also void space within junctions of three or more edges needs to be filled to provide support for following layers.

To meet the basic requirement of a minimum turning radius, there are two approaches: choosing a spline curve or circular arc segments. Strictly parallel curves can be realized with both methods. However, it is important to note that crossings for junctions with more than three edges must be considered and there must not  be an overlap of curves between different connections within these junctions. For this work Bézier spline curves are chosen.

\subsection{Bézier splines}

We utilize both quadratic and cubic Bézier curves. Different mathematical notations and various numerical realizations are established for Bézier curves, depending on the need for numerical stability or performance. We have selected the explicit form, which for $0 \leq t \leq 1$ gives quadratic and cubic Bézier curves as
\begin{align*}
B_2(t) & = (1-t)^2 V + 2\,(1-t)\,t\,U + t^2 W, \\
B_3(t) & = (1-t)^3 V + 3\,(1-t)^2\,t\,S + 3\,(1-t)\,t^2\,T + t^3\,W.\\
\end{align*}
$V,U,W$ and $V,S,T,W$ are control points on a (here) two-dimensional plane, hence $B(t) \in \mathbb{R}^2$. 
The curvature reads as
\begin{align}
\kappa(t) = \frac{ {B^{'}(t)}_x \, {B^{''}(t)}_y - {B^{'}(t)}_y \, {B^{''}(t)}_x }{\left(   {{B^{'}(t)}_x}^2 + {{B^{'}(t)}_y}^2  \right)^{\frac{3}{2}}},
\label{eqn:curvature}
\end{align}
using the $x-$ and $y-$components of the first and second derivatives of the Bézier curves. These are
\begin{align*}
B_2^{'}(t) &= 2\,(1-t)(U-V) + 2\,t\,(W-U), \\
B_2^{''}(t) &= 2\,(W -2\,U + V), \\
B_3^{'}(t) &= 3\,(1-t)^2 (S-V) + 6\,(1-t)\,t\,(T-S) + 3\,t^2(W-T), \\
B_3^{''}(t) &= 6\,(1-t)(T- 2\,S + V) + 6\,t\,(W-2\,T + S).
\end{align*}

\subsection{Quadratic vs. cubic curves}

\begin{figure}[!h]
  \centering
  \includegraphics[width=.45\textwidth]{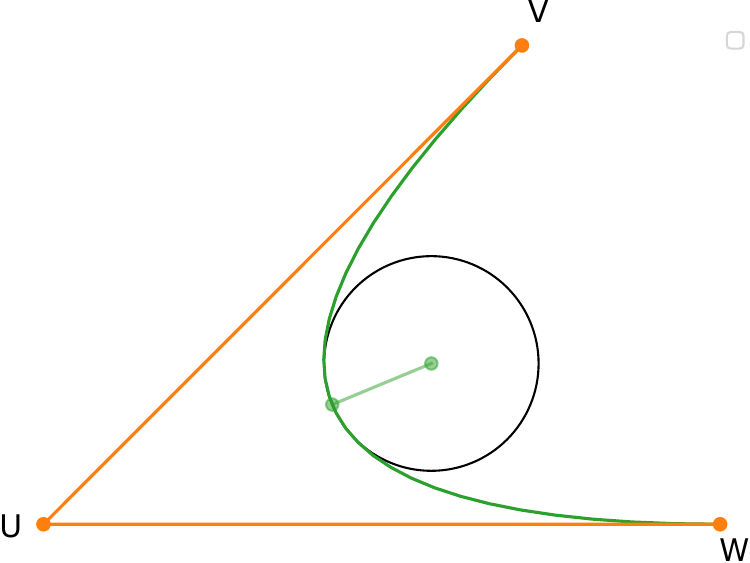} \,
  \includegraphics[width=.45\textwidth]{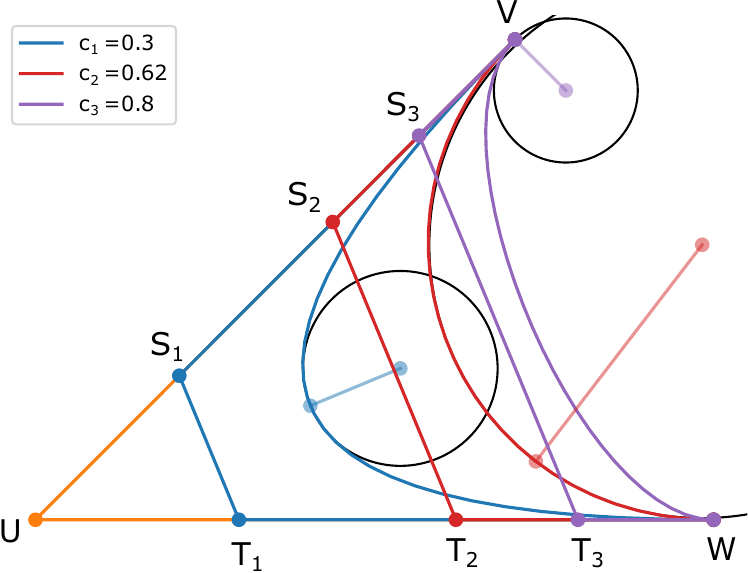} 
  \caption{left: quadratic Bézier curve with three control points $V,U,W$ and an arc segment aligned to the maximal curvature; 
  right: cubic Bézier curves with control points $V,S_i, T_i,W$ with $\overline{U S_i}=\overline{U T_i}=c_i, \;i \in \lbrace 1,2,3\rbrace$. For $c_2$ the red curve is close to a circular arc}
  \label{fig:bezier} 
\end{figure}

\begin{figure}[!h]
  \centering
  \includegraphics[width=.65\textwidth]{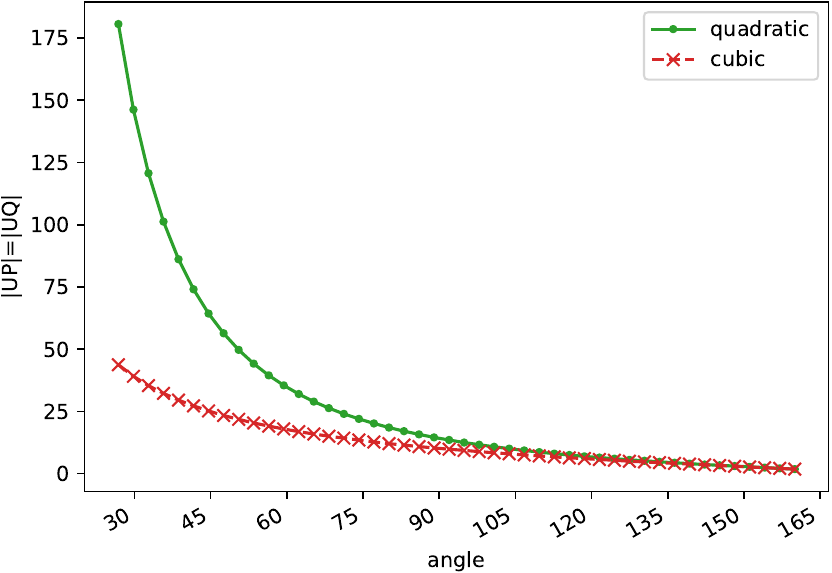}
  \caption{necessary distance $\overline{U V}=\overline{U W}$ in mm to obtain a curve with minimal turning radius 10\,mm with optimal $c$ for the cubic curve.}
  \label{fig:quadratic_vs_cubic} 
\end{figure}

In \figref{fig:bezier} we construct for an isosceles triangle $V,U,W$ with inner angle $45^\circ$ a quadratic Bézier curve. The maximum curvature is indicated by an embedded circle, it has the position $t=0.5$. The maximal curvature can be read as minimum radius for the curve via 
\begin{align*}
R=\frac{1}{\max_t \lvert \kappa(t) \rvert}.
\end{align*}
For a given minimal radius and a given angle, $| U V|$ and $|U W|$ need to be scaled to match the given radius. For   acute angles, quadratic Bézier are inefficient as they require large scaling, see \figref{fig:quadratic_vs_cubic}.

A cubic Bézier curve with its four control points $V,S,T,W$ requires the construction of the inner control points $S$ and $T$. We construct them by using the scaling parameters $c_S$ and $c_T$ 
\begin{align*}
S = U + \frac{c_S}{|U V|} \overrightarrow{U V}, \quad
T = U + \frac{c_T}{|U W|} \overrightarrow{U W}
\end{align*}
For the isosceles triangle we have $c_S=c_T=c$. In \figref{fig:bezier} we see, that there is a $c$ where the curve is very close to a circular arc segment. For smaller $c$, the minimal radius becomes smaller. For larger $c$, parts of the curve get a high curvature at the end regions of $t$. With the requirement of having a curve with minimized maximum curvature, the parameters $c_S$ and $c_T$ are implicitly but indirectly given, and thus the cubic Bézier curve can be constructed. For performance reasons we interpolate an offline computed data set with input data angle and the relationship of the lengths of the line segments $\overline{U V}$ and $\overline{U W}$ and output data $c_V$, $c_W$ and $t$ of the maximal curvature.

In \figref{fig:quadratic_vs_cubic} we show for the isosceles triangle case for different angles the necessary scaling of the line segments to satisfy a 10\,mm turning radius. For sufficient obtuse angles we use quadratic Bézier curves.
  
\begin{figure}[!h]
  \centering
  \includegraphics[width=.45\textwidth]{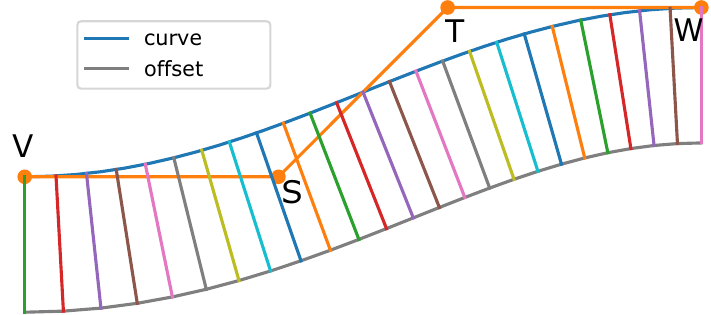}
  \caption{cubic Bézier curve with four control points $V,S,T,W$ and a constructed offset curve via normals.}
  \label{fig:offset} 
\end{figure}

Note that constructing a curve directly parallel to a quadratic or cubic Bézier curve with another Bézier curve is not an easy task. Therefore, it is common practice to discretize the reference curve and create points offset by the distance of the fiber width in the normal direction, using the first derivative of the reference curve, see \figref{fig:offset} and \secref{sec:interlooping}.

\subsection{Constructing Curves}

\begin{figure}[!h]
  \centering
  \includegraphics[width=.5\textwidth]{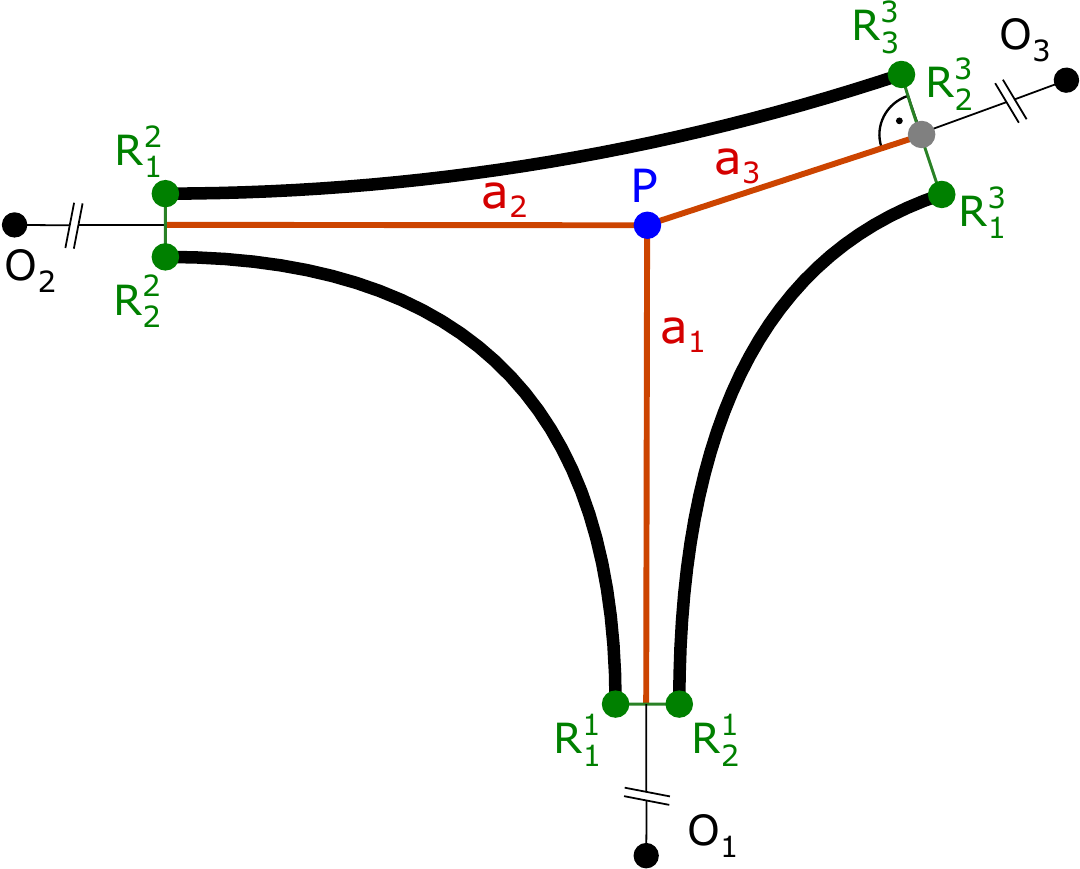}
  \caption{left: Rim points are numbered left to right with a distance $a_i$ from the rim plane to the center $P$.}
  \label{fig:rim} 
\end{figure}

For each junction with its center point $P$, every edge forms a \text{side} where the other endpoint of the edge is denoted as $O$. We denote as \textit{rim} points $R$, where for a fiber path the straight edge connects to the curve. Rim points are arranged on a line segment $\overline{R_1 R_n}$ perpendicular to $\overrightarrow{PO}$ with a distance of $a$ from $P$, see \figref{fig:rim}. The distance $|R_i\,R_{i+1}|$  of the rim points is the printed width of the fiber bundles.  

\begin{figure}[!h]
  \centering
  \includegraphics[width=.45\textwidth]{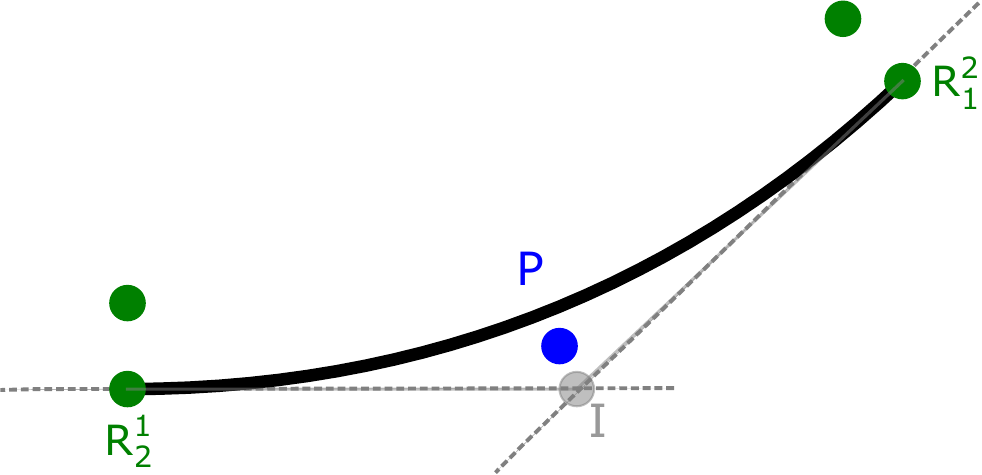} \;
  \includegraphics[width=.45\textwidth]{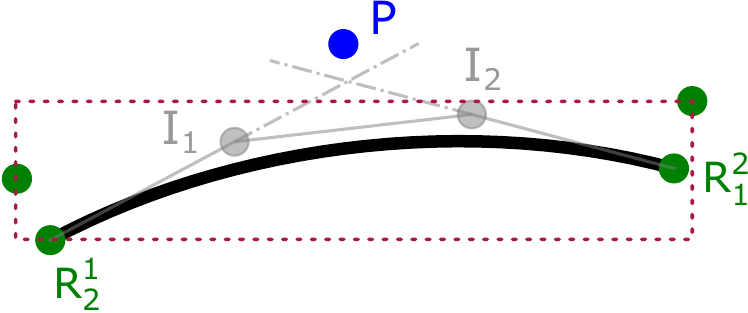} 
  \caption{left: edge bow as quadratic Bézier curve; right: solitary bow as cubic Bézier curve with convex hull depicted as dashed box.} 
  \label{fig:edge_solitary_bow} 
\end{figure}

We refer to the implementation of a connection as \textit{bow}. A bow exists between two rim points, one from each side. We identify three types of bows: edge, solitary and s-shape bow. The Bézier curve parameters for the edge and solitary bows (as shown in \figref{fig:edge_solitary_bow}) are constructed solely from the two constitutive sides. For all bow types, a tangent line is formed for each side using the rim point and $\overrightarrow{PO}$. In case of the \textit{edge bow}, the tangent lines intersection $I$ and the rim points define, depending on the angle, either a quadratic Bézier curve or a cubic curve (the later with implicitly given scaling parameters $c$). If $I$ falls outside the convex hull of the junction (here just the box of minimal and maximal rim point coordinates), a cubic \textit{solitary bow} is constructed.  For the solitary bow, the inner points $I_1$ and $I_2$ are constructed on the tangent lines as
\begin{align*}
  I_1 = R^1 - \frac{\frac{1}{3}}{\lvert R^1 R^2 \rvert}\, \overrightarrow{P O_1}, \quad I_2 = R^2 - \frac{\frac{1}{3}}{\lvert R^1 R^2 \rvert}\, \overrightarrow{P O_2}.
\end{align*}

\begin{figure}[!h]
  \centering
  \subfloat[cubic Beziér curve]{
  \includegraphics[width=.5\textwidth]{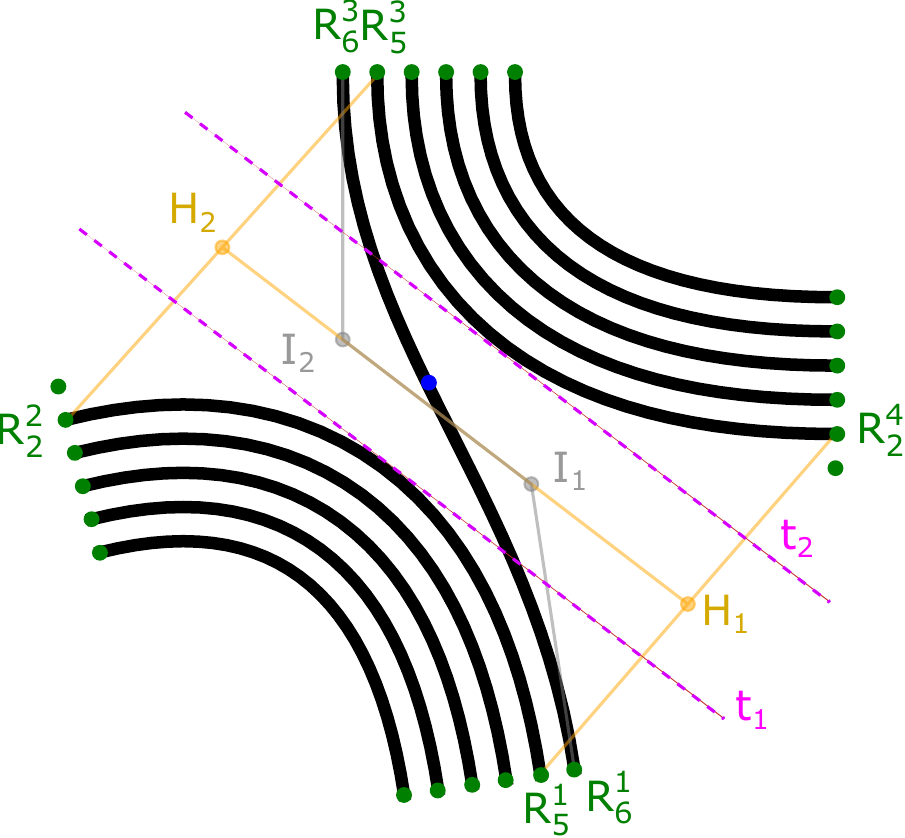}
  \label{fig:s_bow_cubic}} \,
  \subfloat[two quadratic Beziér curves]{
  \includegraphics[width=.45\textwidth]{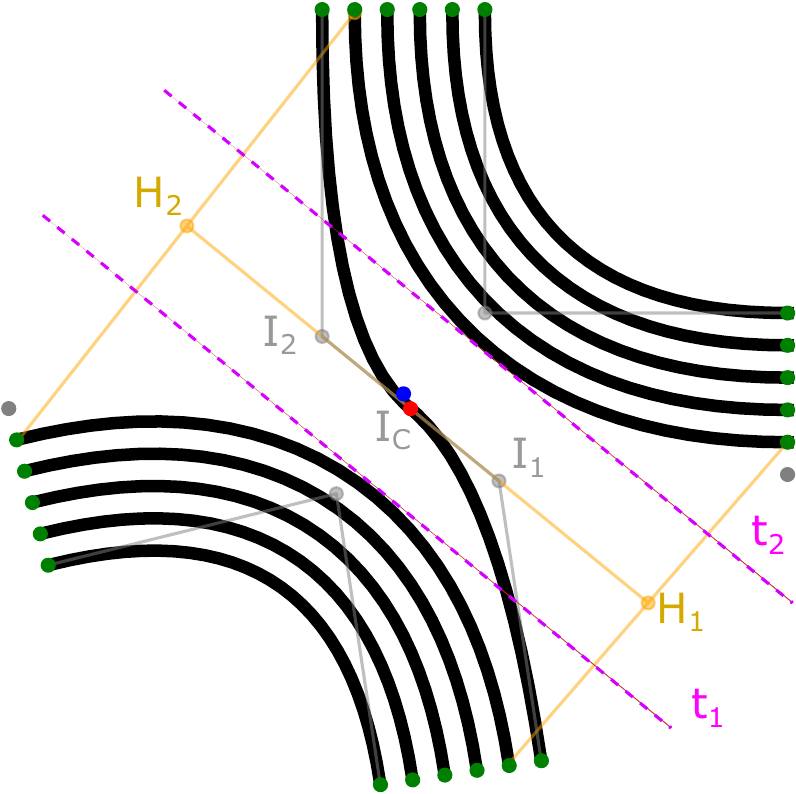} \label{fig:s_bow_biquadratic}}\,
  \caption{An s-shape bow is defined by left and right edge bow. For this we use the cubic Beziér curve shown in (a). An alternative realization via two quadratic Beziér curves, as shown in (b), typically leads to overly large curvature.} 
\end{figure}

The third type of bow is the \textit{s-shape bow} shown in \figref{fig:s_bow_cubic}. It requires two diagonally opposite edge bows (in the example $R_5^1$ to $R_2^2$ and $R_2^4$ to $R_5^3)$. It realizes a crossing connection with at least a fiber bundle width distance to the defining edge bows. 
In essence, a cubic Bézier curve aims to align with the connection of its inner points, namely $I_1$ and $I_2$. The two diagonal opposite edge bows possess a family of parallel tangent lines denoted $t_1$ and $t_2$. 
We approximately select the pair of parallel tangents, that touch the edge bows at their larges curvature. We choose $I_1$ and $I_2$ such, that the line $\overline{I_1 I_2}$ is located parallel and approximately centrally to these tangents. 
This can be achieved by halving the distance of pairs of the edge bow's rim points, resulting in the points $H_1$ and $H_2$. The intersection of the line $\overline{H_1 \, H_2}$ with the tangents of the s-shaped bow's rim points gives the points $I_1$ and $I_2$. 
If there are no distinct $I_1$ and $I_2$ within the convex hull of the junction, the s-shape bow will degenerate into an edge or solitary bow. The idea of aligning the crossing along $\overline{H_1 \, H_2}$ could be more explicitly expressed by introducing a point $I_C$ between $I_1$ and $I_2$ and using two quadratic Beziér curves for the points $R_6^1, I_1, I_C$ and $I_C, I_2, R_6^3$, respectively (see \figref{fig:s_bow_biquadratic}). However, this exhibits in many cases too large curvatures.

In the view of this paper and the reference implementation, we assume that junctions have no more than four edges. With proper sheet definitions there can be at most one s-shape bow at junction (within a layer). For each connection type we construct a reference bow. First, edges use the first and last rim points for reference bows and then fill up with offset edge bows. Degenerated edge bows are realized via solitary bows.  Next, center crossing bows can be constructed as reference bow with offset bows. A crossing without diagonal edge bows is realized as solitary bow.
    
\subsection{Parameterizing junctions}
In this section we describe our approach to find for junctions the proper bow parameters, such that all paths comply with the minimal turning radius. What needs to be determined is the distance of the rim points to the junction's center point $P$. Each edge of a junction forms a side, which contains the rim points. Two adjacent sides form a \textit{wedge}. Each wedge has a left rim point $R^\text{l}$ and a right rim point $R^\text{r}$. The side's tangents through the rim points intersect for each wedge in a point $U$. Each wedge has an angle $\beta = \angle\,R^\text{l} U R^\text{r}=  \angle \,O^\text{l} P O^\text{r}$.

\begin{figure}[!h]
  \centering
  \includegraphics[width=.4\textwidth]{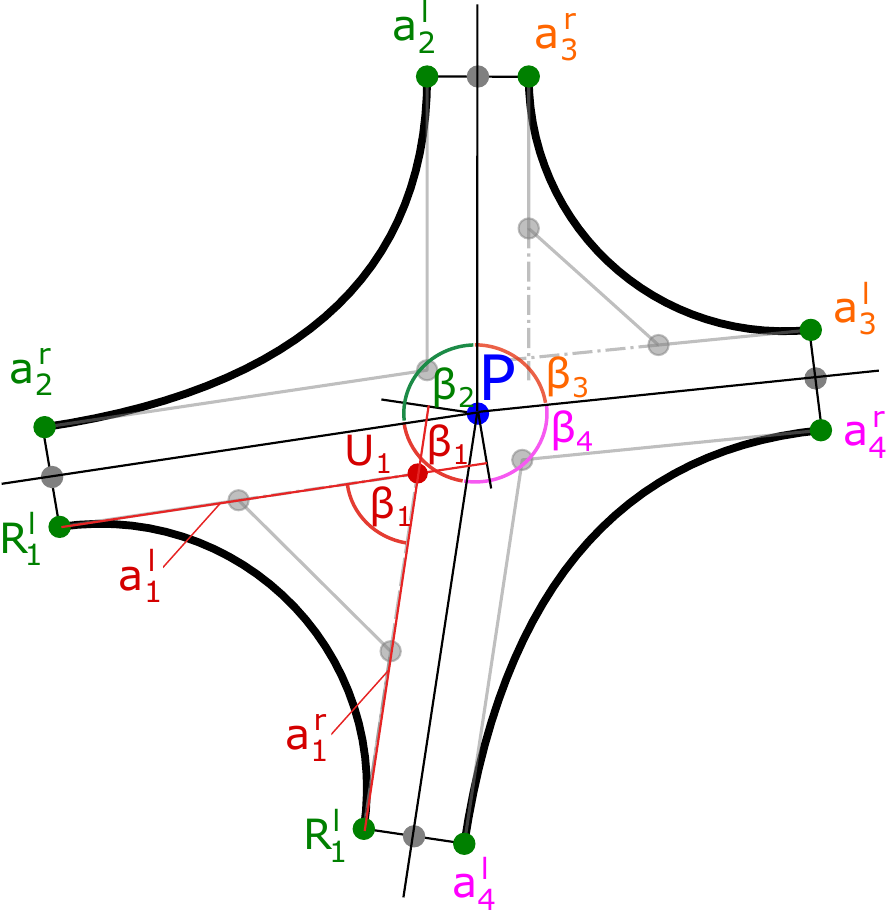} 
  \caption{Determining the $a$-distance values. Each of the four wedges consists of a left and right $a$ distance and an angle.} 
  \label{fig:angles} 
\end{figure}

In \figref{fig:rim} we denote the  distance of the junction's center point $P$ to the rim point by $a$. We can write this value more detailed as $a = v + u$, with $v=|R\, U|$ and $u$ is the distance of $U$ projected to $O P$ and $P$. Each wedge has one point $U$ but a left and right distance $u^\text{l}$ and $u^\text{r}$ which might be different for different edge widths.      
  
Our procedure for junctions of three and four edges is the following: First we assume for each wedge $v^\text{l}=v^\text{r} =: v$ and find the $v$ which leads to an edge bow with the required minimal turning radius. Depending on the angle this is either a quadratic or cubic Bézier curve. Then we find the wedge which has the largest $a$ distance and fix this wedge. Commonly this is the wedge with the smallest angle, but with significantly different edge widths, this might not be the case. In the example this is wedge $1$ and we obtain $a_1^\text{l}=v_1 + u_1^\text{l}$ and $a_1^\text{r}=v_1 + u_1^\text{r}$. 

$a^\text{l}$ and $a^\text{r}$ of the fixed wedge are then fixed for the counterparts of the adjacent wedges leaving a single $v^\text{r}$ and $v^\text{l}$ of the two adjacent wedges to be determined to satisfy the minimal turning radius. In the example this fixes $a_2^\text{r} := a_1^\text{l}$ and $a_4^\text{l} := a_1^\text{r}$ and we search the determined $v_2^\text{l}$ and $v_4^\text{r}$  
  
In the case of a junction of four edges, we have a candidate for the remaining wedge 3 with $a_3^\text{r}=\gamma \,v_3^\text{r}+u_3^\text{r}$ and $a_3^\text{l}=\gamma \,v_3^\text{l}+u_3^\text{l}$. We then find the first $\gamma >= 1$ such that the minimal turning radius is not violated for the remaining wedge. 

Finally, we check if adjacent $a$ distances differ and correct the smaller value in case.

The described algorithm consists of a sequence of minimization problems, where only a single parameters is to be found for each problem (e.g. via bisection). This makes the procedure computationally quite efficient.

We ensure, that the curvature is not violated for the outer reference edge curves. All offset edge curves are parallel and have smaller curvature. However, we do not check the construction of the reference solitary and s-shaped bows and their offsets. We also have no proof for the critical requirement to ensure fiber paths do not overlap with other bows. No issues were observed in our numerical experiments, as long as the minimal radius is sufficiently large. In doubt, one could verify the conditions numerically and adjust the curvature parameter accordingly for each junction. 

The only parameters for the fiber path planning algorithm are the minimum turning radius and fiber width. All remaining variables are determined by the fiber path planning algorithm to meet the required curvature and maintain a good parallelism between fiber paths.

\subsection{Interlooping}
\begin{figure}[!h]
  \centering
  \includegraphics[width=.45\textwidth]{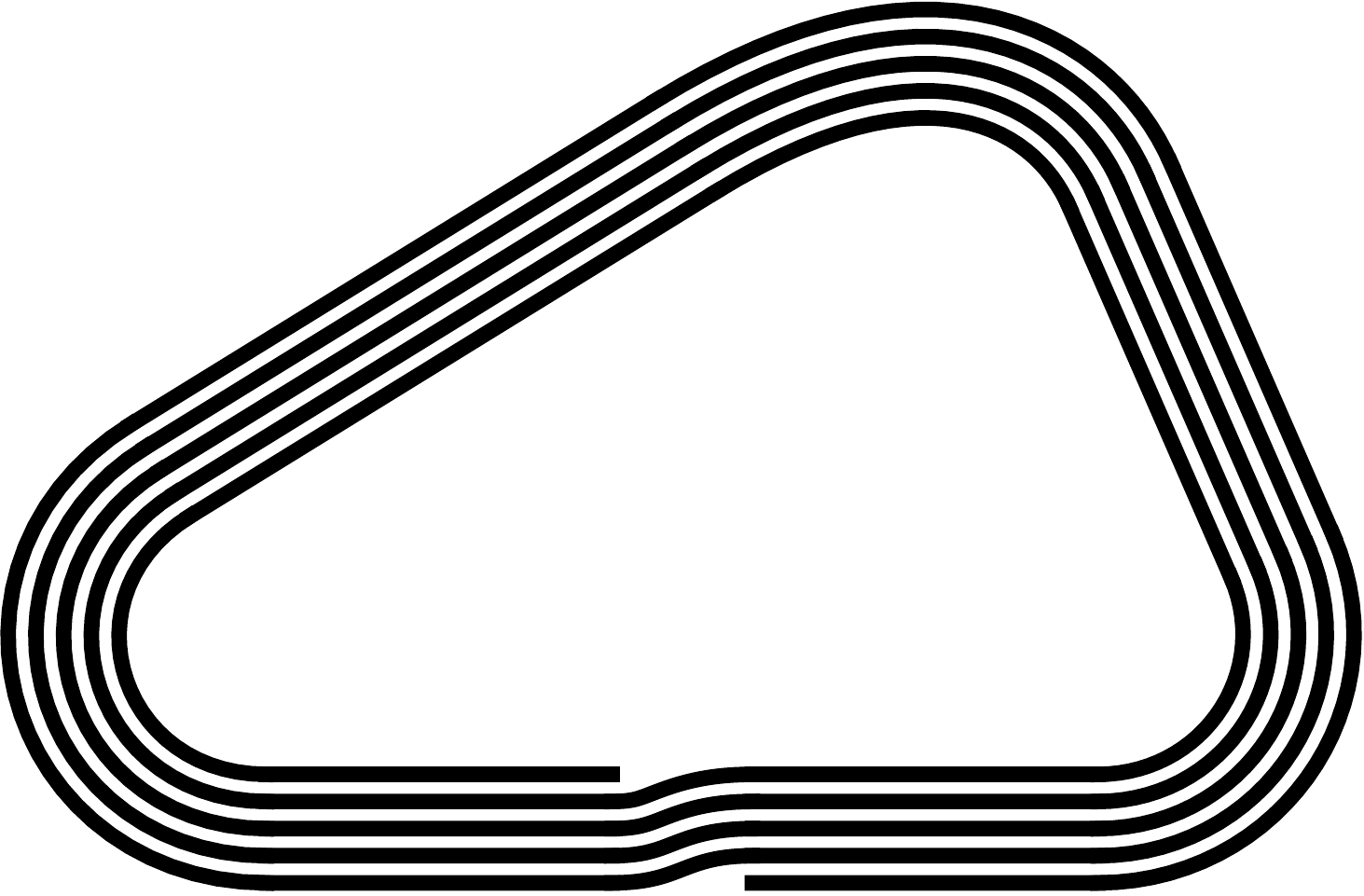}
  \caption{\label{fig:interlooping}Interlooping: connecting instances of a concentric closed loop to a single path. Each interconnecting curve is different, see \figref{fig:interloop_curves}.}
\end{figure}

We differentiate between open and closed loops. The fiber bundles of open loops have two end point coordinates for the printing process. %
For concentric closed loop instances, we suggest to connect the fiber bundles  to a single continuous path, see \figref{fig:interlooping}. We call this \textit{interlooping} and obtain a single pair of end points for the loop set.

\begin{figure}[!h]
  \centering
  \subfloat[curves]{
  \includegraphics[width=.4\textwidth]{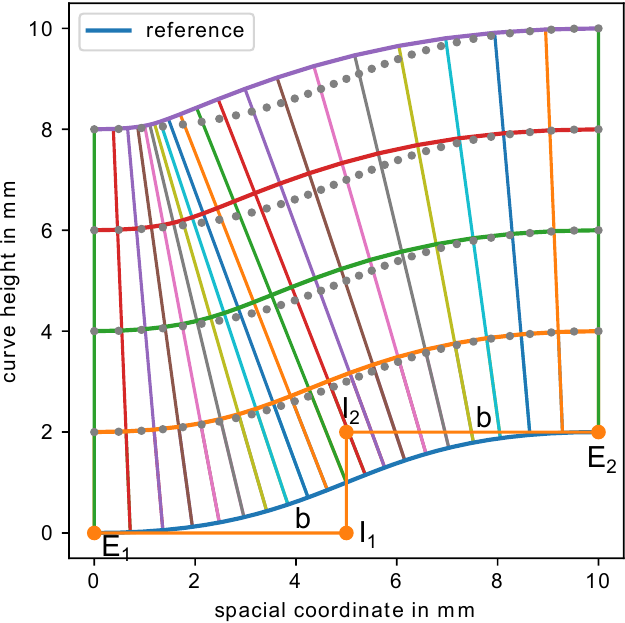}
  \label{fig:interloop_curves} }\,
  \subfloat[curvature]{
  \includegraphics[width=.5\textwidth]{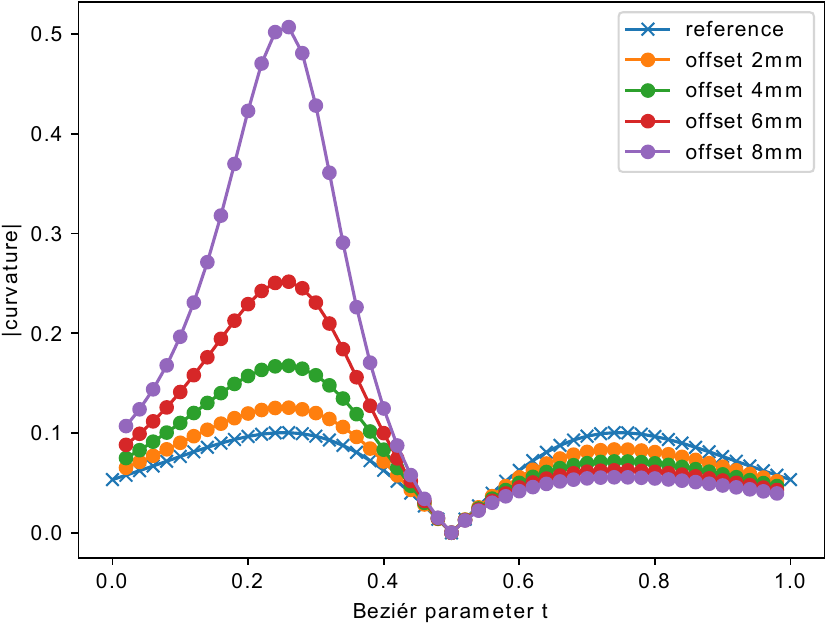} 
  \label{fig:interloop_curvature} }
  \caption{For interlooping, we show the cubic Beziér reference curve and 4 offset curves assuming 2\,mm fiber width and 10\,mm minimal turning radius for the reference curve. In (a) we show constructed Beziér curves by dotted lines, however these are not parallel.}
\end{figure}

To respect the minimal turning radius for interlooping, we construct a cubic Beziér curve from the points $E_1, I_1, I_2$ and $E_2$ as depicted in \figref{fig:interloop_curves}. The control polygon has $90^\circ$ angles at $I_1$ and $I_2$, the distance $|I_1 \, I_2|$ is the fiber width (here 2.0\,mm). We introduce a single parameter $b$ for the distance $|E_1\,I_1|$ and $|I_2\,E_2|$. Numerically, we can easily find $b$ to respect the minimal turning radius (here 10\,mm) as $b=5.0$. To print parallel interlooping curves, we need offset curves, as no parallel cubic Beziér curves can be constructed (see the dotted constructed cubic curves in \figref{fig:interloop_curves}). However, we need to handle the increasing curvature of offset curves. While the curvature for a cubic Beziér curve is  directly given by \eqnref{eqn:curvature}, we obtain it for offset curves by the radius of the circle defined by three adjacent offset points. We plot in \figref{fig:interloop_curvature} the curvature and not the reciprocal value radius, as the curvature is zero in between $I_1$ and $I_2$.

We numerically search $b$ such that the minimal radius is satisfied for the largest offset (can be efficiently performed using bisection). For offsets of 1, 2, 3, 4 and 5 fiber bundles, $b$ is 5.5, 6.0, 6.4, 6.9 and 7.2\,mm. One can also locate the reference curve in between and have positive and negative offsets.

\subsection{Loop identification}
\label{sec:interlooping}
\begin{figure}[!h]
  \centering
  \includegraphics[width=.35\textwidth]{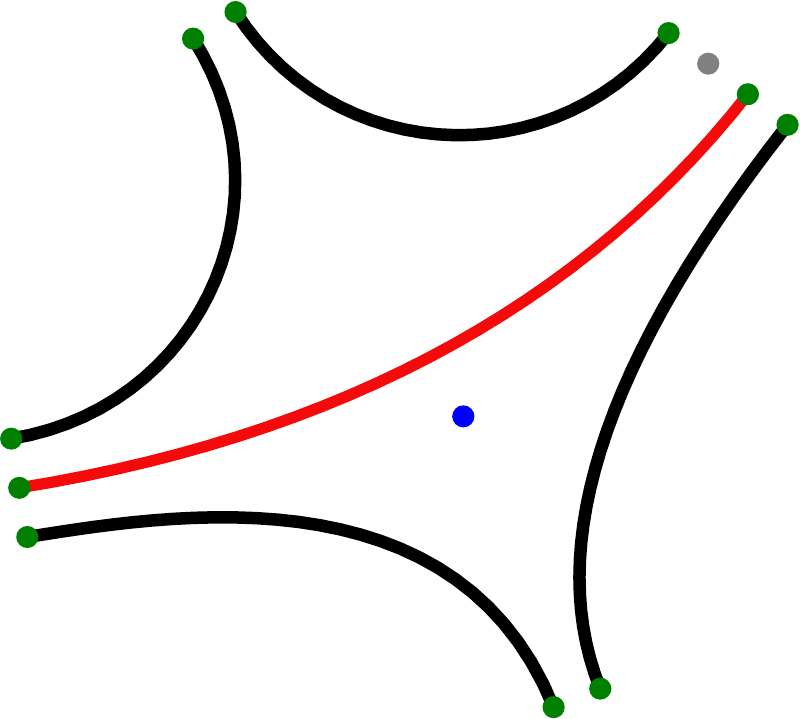} 
  \caption{The right position of the red crossing is ambiguous. Knowledge of loop assignments in the adjacent junctions is required to resolve this.} 
  \label{fig:ambiguous} 
\end{figure}

In principle, it would not be necessary to track the relationship of connections to loops in the path planning implementation, but in certain configurations there is a need. An example is given in \figref{fig:ambiguous}, where the position of the crossing is ambiguous. To resolve it, the positions of the corresponding loop in the other junctions needs to be known. The solution is rather technical and considered there as implementation detail but not part of the description of the principal algorithmic approach itself.

\section{Numerical example}
\label{sec:example}

\begin{figure}[!h]
  \centering
  \subfloat[topology optimization]{
  \includegraphics[width=.45\textwidth]{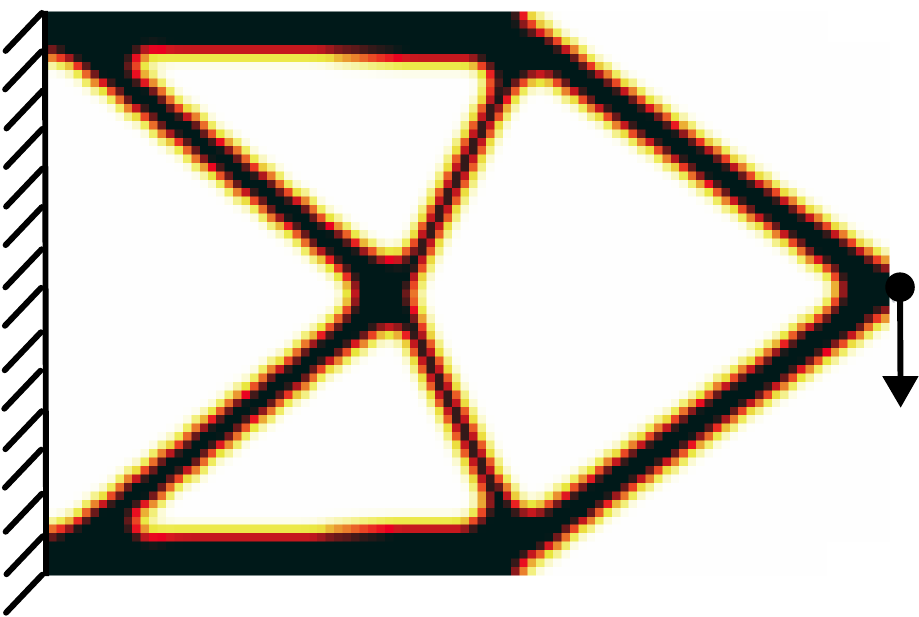} \label{fig:itop}} \;
  \subfloat[feature-mapping]{
  \includegraphics[width=.45\textwidth]{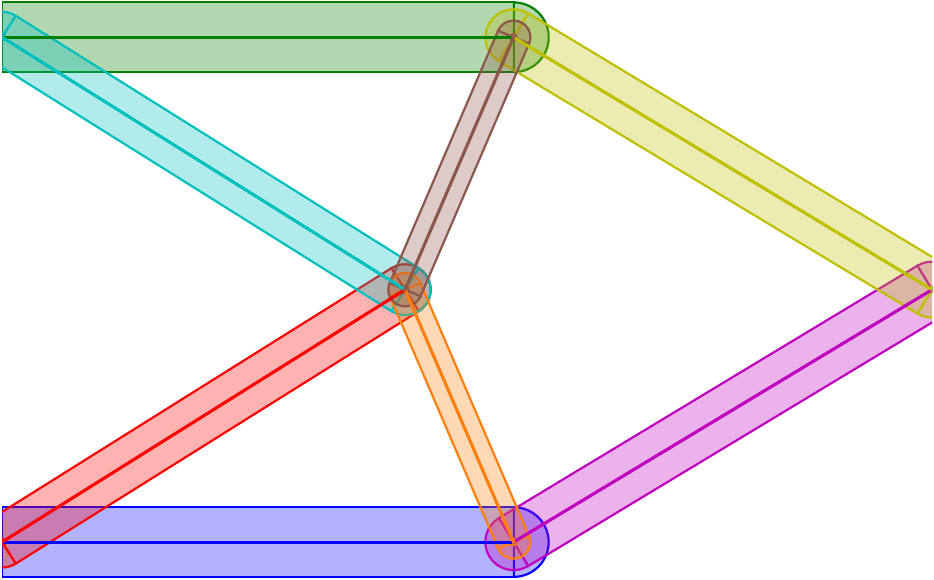} \label{fig:opt_profiles} }
  \captionof{figure}{We perform topology optimization for a standard one-load cantilever problem. Via feature-mapping the bar widths are optimized.}
  \label{fig:cantilever_setup} 
\end{figure}

Density based topology optimization (a.k.a. SIMP optimization, see \cite{Bendsoe:2003:Book}) is a standard approach to perform structural optimization. Several CFRP publications use numerical examples based on designs resulting from topology optimization. A classical topology optimization benchmark problem is the cantilever shown in \figref{fig:itop}. The design was obtained by the web application iTop (\url{monster.mi.uni-erlangen.de/iTop}) but other academic tools, e.g. from \url{www.topopt.mek.dtu.dk/apps-and-software}, or several finite element package can obtain similar designs.

For the layered fiber patterns, we need a graph representing the design. Truss optimization is an another topology optimization approach, which directly results in graphs (see again \cite{Bendsoe:2003:Book}). Here, we chose with feature-mapping a rather modern structural optimization approach, which to some extent combines the two mentioned ones. The design variables directly describe bar elements, but the continuum based finite element analysis is more accurate in resolving the connection of bars compared to the reduced truss model. For an overview of feature-mapping, we refer to the review paper \cite{Wein:2020:Review}. In \figref{fig:opt_profiles} we show a manual interpretation of the design in \figref{fig:itop} by bar elements. For simplicity, endpoints of the bars are fixed, but the thickness parameter of the bars is optimized. We use the feature-mapping approach presented in \cite{Greifenstein:2023:Spline} via the reference implementation in openCFS (\url{gitlab.com/openCFS}) but assume for isotropic material.

\begin{figure}[!h]
  \centering
  \includegraphics[width=.35\textwidth]{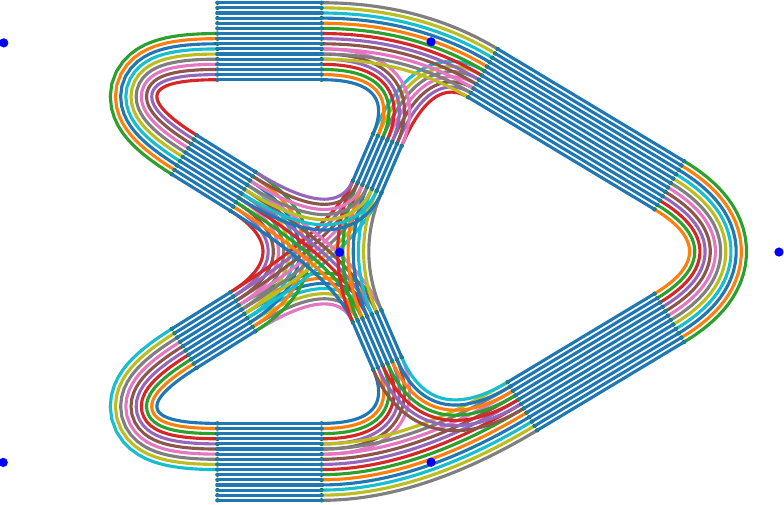} \;
  \includegraphics[width=.35\textwidth]{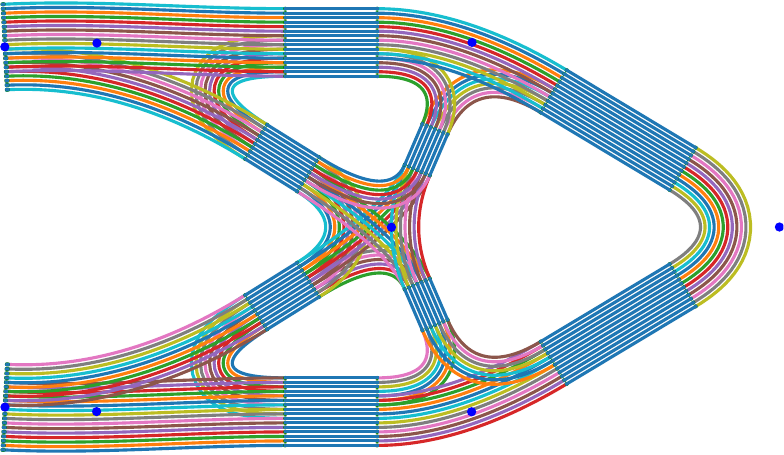} 
  \caption{Set of possible connections for the graph in \figref{fig:opt_profiles} without and with support extension. Note the support nodes on the left and load point on the right.} 
  \label{fig:cantilever_struct} 
\end{figure}

Applying the path planning, we see in the left figure in \figref{fig:cantilever_struct}, that the left support points and the right load point are not reached by the fibers. For the support points with the rather acute angles this is even more prominent. Furthermore, the horizontal and diagonal edge at the support points have no matching width. The latter is due to the fact that support points in structural optimization are not just points, where edges connect, but the edges are intended to connect with the support (fixed displacement). We therefore add further one-node junctions (stubs) left of the support points to our graph. The width of this stub extensions is chosen to be a little more than the width of the horizontal edges.  
One could question the the necessity of the edge bows connecting edges of support points. 

While layered fiber pattern optimization can easily be tuned by modifying connection targets $\tilde{c}$, this would ask for a numerical validation of the settings, which is outside the scope of this work. We do not need an extension for the load point, as horizontal oriented fiber bundles for a vertical load does not make  much sense.       

\begin{figure}[!h]
  \centering
  \includegraphics[width=.45\textwidth]{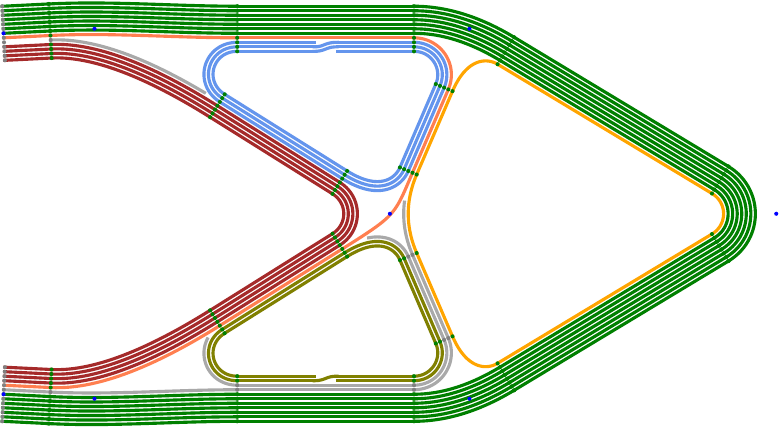} \;
  \includegraphics[width=.45\textwidth]{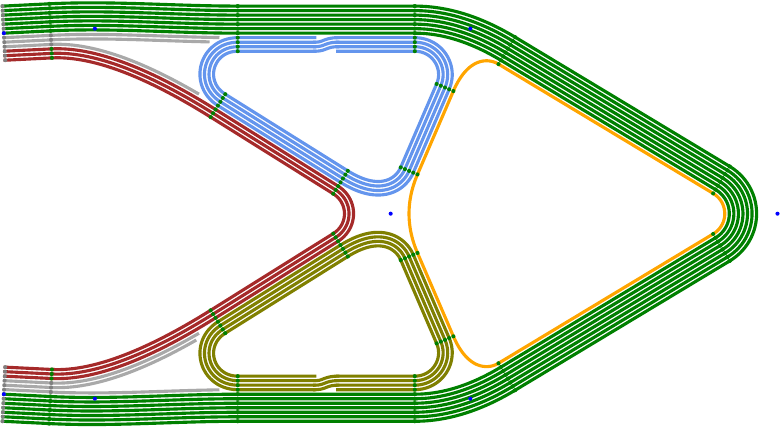} 
  \caption{Obtained layers for the cantilever problem. The left layer (first layer) occurs also vertically flipped with another sheet. We also visualize the vertices and rim points. In gray we fill unused loop space. Interlooping connects instances of closed loops.} 
  \label{fig:cantilever_layers} 
\end{figure}

We optimize the layered fiber patterns with penalty $p=2$ and the constraint to use all edge bows \eqnref{eqn:Cs_x_gt_cb}. Three different types of patterns are obtained, with two of them mirrored variants of the same pattern, see \figref{fig:cantilever_layers}. Without enforced edges, a total of 6 different patterns are obtained for the first 100 layers, where three patterns are dominant. For a structure with more volume and as such more fiber bundles (see \figref{fig:cantilever_struct}) we still obtain six different patterns but much more uniformly distributed.

\section{Discussion and possible extensions}
\label{sec:discussion}
In this paper, we present a novel approach to generate paths for continuous fiber-reinforced fiber bundles. The pattern of these paths is found by solving multiple discrete linear programming problems. These problems must be parameterized, but this parameterization makes it possible to map even complex structural requirements onto the patterns. 

Key components of our approach are the introduction of:
\begin{enumerate}[label=(\roman*)]
\item The use of discrete optimization problems to find optimal fiber patterns.
\item Taking advantage of the fact that multiple layers are typically printed to satisfy complex and conflicting requirements.
\item Generate spline based fiber paths that satisfy a minimum turning radius and are strictly parallel.
\end{enumerate}

The aforementioned pattern optimization parameterization will be published in future work. Initial tests have already shown that connection targets $\tilde{c}$ can be found using an FE-based optimization problem. We also plan to develop an algorithm to automatically determine the relevant loops.

The optimization of layered fiber patterns is inherently three-dimensional, but the proposed path planning is implemented for two-dimensional layers only. Pattern optimization and path planning can also be applied individually by substituting the counterpart by other approaches.  

We added the option to optionally fill void edges and full and partial bows. fillings are outside loop configurations and usually do not have a continuous connection from load to support points. It is up to the user to decide between fill some space by C-CFRTP or standard isotropic filling material.   

\section*{Replication of results}
\noindent We provide open-source (MIT/X11 license) reference implementations in Python for the discrete layered fiber pattern optimization and the spline based fiber path planning at \texttt{gitlab.com/fiona3644621/FionaPathLayout}.

\section*{Acknowledgements} 
\noindent The authors gratefully acknowledge the financial support by the German Federal Ministry for Economic Affairs and Climate Actions (BMWK) in the course of the FIONA (LuFo IV-1, FKZ: 20W1913F) project.\\

\bigskip
\noindent During the preparation of this work the authors used deepl.com in order to improve the language. After using this tool/service, the authors reviewed and edited the content as needed and take full responsibility for the content of the publication.

\bibliographystyle{plainnat}      %
\bibliography{optimization, cfrp} 

\end{document}

%% file: fabian_macros.tex
\usepackage{amsmath,amssymb,amsfonts,amsbsy}
\usepackage{mathptmx}
\usepackage{setspace}
\usepackage{bm}

\newcommand{\tabref}[1]{Tab.~\ref{#1}}
\newcommand{\figref}[1]{Fig.~\ref{#1}}
\newcommand{\secref}[1]{Sec.~\ref{#1}}

\newcommand{\eqnref}[1]{(\ref{#1})}